\documentclass[twoside,11pt]{amsart}
\usepackage{latexsym}
\usepackage{amssymb}

\makeatletter
\def\serieslogo@{}
\makeatother \makeatletter
\def\@setcopyright{}
\makeatother

\def\ep{\epsilon}

\newtheorem{thm}{Theorem}[section]

\newtheorem{lem}[thm]{Lemma}

\theoremstyle{definition}

\theoremstyle{remark}
\newtheorem{rem}{Remark}[section]

\numberwithin{equation}{section}

\title[Recovery of high frequency wave fields]{Recovery of high frequency wave fields from phase space based measurements}

\author{Hailiang Liu and James Ralston}
\address{Iowa State University, Mathematics Department, Ames, IA 50011} \email{hliu@iastate.edu}
\address{UCLA, Mathematics Department, Los Angeles, CA 90095} \email{ralston@math.ucla.edu}


\date{May 23, 2009}
\begin{document}

\begin{abstract}Computation of high frequency solutions to wave equations is important in
many applications, and notoriously difficult in
resolving wave oscillations. Gaussian beams are
asymptotically valid high frequency solutions
concentrated on a single curve through the
physical domain, and superposition of Gaussian
beams provides a powerful tool to generate more
general high frequency solutions to PDEs. An
alternative way to compute Gaussian beam
components such as phase, amplitude and Hessian
of the phase, is to capture them in phase space
by solving Liouville type equations on uniform
grids. In this work we review and extend recent
constructions of asymptotic high frequency wave
fields from computations in phase space. We give
a new level set method of computing the Hessian
and higher derivatives of the phase. Moreover, we
prove that the $k^{th}$ order phase space based
Gaussian beam superposition converges to the
original wave field in $L^2$ at the rate of
$\ep^{\frac{k}{2}-\frac{n}{4}}$ in dimension $n$.

\end{abstract}

\maketitle

\bigskip

\tableofcontents

\section{Introduction}
In this paper we consider the following  equation
\begin{align}\label{pp}
P\psi =0, \quad (x, t)\in \mathbb{R}^n \times \mathbb{R},
\end{align}
where  $P=-i\epsilon \partial_t +H(x, -i\epsilon
\partial_x)$ is a linear differential operator with a real principal symbol $\tau + H(x, p)$,
subject to the highly oscillatory initial data
\begin{equation}\label{ini}
\psi(x, 0)=\psi_{in}(x):=A_{in}(x) e^{iS_{in}(x)/\ep},
\end{equation}
where $A_{in} \in C_0^\infty(\mathbb{R}^n)$ and
$S_{in}\in C^\infty(\mathbb{R}^n)$. The canonical
example is the semi-classical Schr\"{o}dinger
equation with the Hamiltonian $ H(x,
p)=\frac{1}{2}|p|^2+V_{ext}(x)$, where
$V_{ext}(x)$ is a given external potential. The
small parameter $\ep$ represents the fast space
and time scale introduced in the equation, as
well as the typical wave length of oscillations
of the initial data. Propagation of oscillations
of wave length $O(\ep)$ causes mathematical and
numerical challenges in solving the problem. In
this article we are interested in the
construction of globally valid asymptotic wave
fields and the analysis of their convergence to
the true  solutions of the initial value problem.

Geometric optics, also known as the WKB method or
ray-tracing, when applied to model high frequency
wave propagation problems such as (\ref{pp})
leads to the WKB-type system for both phase and
amplitude. The phase is governed by the
Hamilton-Jacobi equation
\begin{equation}\label{ss}
   \partial_t S +H(x, \nabla_x S)=0, \quad x\in \mathbb{R}^n, \; t>0.
\end{equation}
Solving this equation using the method of
characteristics can lead to singularities which
invalidate the approximation. In general, this
breakdown occurs when the density of rays becomes
infinite. This corresponds to the formation of a
caustic where geometric optics incorrectly
predicts that the amplitude of the solution is
infinite. The consideration of these
difficulties, beginning with Keller
\cite{Keller:1958}, and Maslov
\cite{Maslov-Fedoriuk:1981}, led to the
development of the theory of Fourier integral
operators, e.g., as given by H\"{o}rmander
\cite{HOb:1971}.

A closely related alternative to Fourier integral
operators is the construction of approximations
based on Gaussian beams. Gaussian beams are
asymptotic solutions concentrated on classical
trajectories for the Hamiltonian $H(x,p)$, and
they have a history going back to at least the
late 1960's. They were initially used to study
resonances in lasers \cite{BB:1972}, and later to
obtain results on the propagation of
singularities in solutions of PDE's
\cite{HOa:1971} and \cite{Ralston82}. At present
there is a considerable interest in using
superpositions of beams to resolve high frequency
waves near caustics. This goes back to the
geophysical applications in \cite{CPP:1982} and
\cite{Hill01}. Recent works in this direction
include \cite{TQR:2007} on \lq\lq gravity waves",
\cite{LQ:2008} and \cite{JWY:2008} on the
semi-classical Schr\"odinger equation, and
\cite{Tan08} and \cite{MR:2008} on acoustic wave
equations.

An alternative to the standard WKB method is to
use multi-valued solutions, $\{S_i(t,
x)\}_{i=1}^K$, to (\ref{ss}) corresponding to
crossing waves \cite{SMM:2003}. This is in sharp
contrast to the notion of viscosity solution
\cite{CL:1983} adopted when (\ref{ss}) arises in
other applications. In the last decade a
considerable amount of work has been done to
 capture multi-valued phases
associated to the WKB system numerically; we
refer to review articles
 \cite{ER:2003, Ru:2007} and references therein. Recently the level set method has been
  developed to resolve
 multi-valuedness of involved quantities in phase space as well as to compute physical observables, e.g.
  \cite{OCKST:2002, CLO:2003, JO:2003,QL:2004,JLOTa:2005, JLOTb:2005,
  LCO:2006, LOT:2006, LWb:2007, LW:2008,
  LW:2009}. The key idea,  for instance in \cite{ CLO:2003, LCO:2006}, is to represent characteristic
 trajectories by common zero sets of some implicit level set functions, and evolve all relevant
 quantities in phase space, see the review \cite{LOT:2006}.
Another phase space based approach is the use of
the Wigner transformation \cite{Wigner:1932} to
map solutions of the underlying wave equation to
functions on phase space. See \cite{LP:1993,
MMP:1994, MPP:1997,SMM:2003} for the application
of the Wigner transformation to the
semi-classical Schr\"{o}dinger equation. These
phase space based approaches are extremely useful
since they unfold `caustics'. However, at
caustics, neither gives correct prediction for
the amplitude. This has led to an effort to
combine the accuracy of beam superpositions at
caustics with the level set method of computation
in phase space. This approach with
attention to the accuracy of the resulting
approximations is the subject of this paper.

In this article we have two objectives:
\begin{itemize}
\item[(i)] to present the construction of beam superpositions by level set methods
;
         \item[(ii)]to estimate the error between the exact wave field and the asymptotic ones.
                                                                                 \end{itemize}
The construction for (i) is based on Gaussian
beams in physical space, but it is carried out by
solving inhomogeneous Liouville equations in
phase space as in \cite{LQB:2007}, \cite{LQ:2008}
and \cite{JWY:2008}. The result is no longer a
superposition of asymptotic solutions to the wave
equation (1.1)! Since it can be written as a
superposition of standard Gaussian beams composed
with a time-dependent symplectic change of
variables, the {\it superposition over phase
space } is still an asymptotic solution, as was
pointed out in \cite{LQB:2007}. Here we consider
{\it superpositions over subdomains} moving with the
Hamiltonian flow, and show directly that they are
asymptotic solutions without reference to
standard Gaussian beams. For (ii) we use the
well-posedness theory for (1.1), i.e. the
continuous dependence of solutions of $P\psi=f$
on their initial data and $f$. Thus, the sources
of error in the Gaussian beam superposition for
an initial value problem are the error in
approximating the initial data and the error in
solving the PDE.

 To be specific, our asymptotic solution is expressed as
                                                                  \begin{equation}\label{sol1}
\psi^\epsilon(t, y)=Z(n,\epsilon)\int_{\Omega(t)}\psi_{PGB}(t,
y,X)dX,
\end{equation}
where $X=(x, p)$ denotes variables in phase space
$\mathbb{R}^{2n}$, $\Omega(0)$ is the domain
where we initialize  Gaussian beams from the
given data, and $\Omega(t)=X(t, \Omega(0))$ is
the image of $\Omega(0)$ under the Hamiltonian
flow. Here $\psi_{PGB}(t,y,X)$ is the phase space
based Gaussian beam {\it Ansatz}, and $Z(n ,\ep)$
is a normalization parameter chosen to match
initial data against the Gaussian profile.  Our
result shows that for the $k^{th}$ order phase
space Gaussian beam superposition, the following
estimate holds on any bounded time interval,
$|t|\leq T$,
\begin{equation}\label{est1}
\|(\psi^\epsilon -\psi)(t, \cdot)\|_{L^2} \lesssim\|(\psi^\epsilon(0, \cdot)-\psi_{\rm in}(\cdot) \|_{L^2} + |\Omega(0)|\ep^{\frac{k}{2}-\frac{n}{4}}.
\end{equation}
Here and in what follows we use $A \lesssim B$ to denote the estimate $A\leq CB$ for a constant $C$ which is
independent of  $\ep$. 

For the initial data of the form $\psi_{\rm
in}=A_{\rm in}(x) e^{iS_{\rm in}(x)/\ep}$ we need
a superposition over an n-dimensional submanifold
of phase space. The asymptotic solution is then
represented as
\begin{equation}\label{sol2}
\psi^\epsilon(t, y)=Z(n,\epsilon)\int_{\Omega(t)}\psi_{PGB}(t,
y,X)\delta(w(t, X))dX,
\end{equation}
where $w$ is obtained from the Liouville equation
$$
\partial_t w+H_p\cdot \nabla_x w-H_x \cdot \nabla_p w=0, \quad w(0, X)=p-\nabla_x S_{\rm in}(x).
$$
Our result shows that
\begin{equation}\label{est2}
\|(\psi^\epsilon -\psi)(t, \cdot)\|_{L^2} \lesssim |{\rm supp}(A_{\rm in})|\ep^{\frac{k}{2}-\frac{n}{4}}.
\end{equation}

We now conclude this section by outlining the rest of this paper: in Section 2 we start with
Gaussian beam solutions in physical space, and
define the phase space based GB ansatz through
the Hamiltonian map. Section 3 is devoted to a
recovery scheme through superpositions over a
moving domain.  The total error is shown bounded
by an initial error and the evolution error of
order $\ep^{1/2-n/4}$.  Control of initial error
is discussed in Section 4, followed by the
convergence rate obtained for first order GB
solutions. In Section 5, we discuss how to use
caustic structure to obtain some better error
estimates. Both convergence and convergence rate
are obtained for higher order GB solutions in
Section 6. In Section 7, we present a new level
set approach for construction of the phase and
their derivatives. Finally in the appendix, we
derive phase space equations for all involved
Gaussian beam components.

\section{Phase space based Gaussian beam Ansatz}

\subsection{First order Gaussian beam solutions}
As is well known,  the idea underlying Gaussian beams (GB) \cite{Ralston05} is to build asymptotic
solutions concentrated on a single curve in physical space $(t, x)\in \mathbb{R}\times \mathbb{R}^n$.
This means that, given a curve $\gamma$ parameterized by $x = x(t)$, one
makes the Ansatz
\begin{equation}\label{an}
\psi^\epsilon(t, x)=A^\epsilon(t, x)e^{i\Phi(t, x)/\epsilon},
\end{equation}
where $\Phi(t, x(t))$ is real, and $Im\{\Phi(t, x)\} > 0$ for $x \not=
x(t)$. The amplitude is allowed to be complex and have an asymptotic expansion in terms of $\ep$:
$$
A^\ep(t, x)=A_0(t, x) +\ep A_1(t, x)+\cdots+\ep^N A_N(t, x).
$$
We wish to build
asymptotic solutions to $P\psi(t, x) = 0$, i.e., we want $P\psi^\epsilon = O(\epsilon^{2})$. Substituting from (\ref{an}),
$$
 P \psi^\ep  = e^{i\Phi(t, x)/\epsilon} \left[ (\partial_t \Phi +H(x, \partial_x \Phi))A_0 + \ep(-iLA_0+(\partial_t \Phi +H(x, \partial_x \Phi))A_1) \right]+O(\ep^2),
 $$
where $L$ is a linear differential operator, whose form is clear from (\ref{ao}) below.  The key step in the GB construction is the choice of $\Phi$ such that $\partial_t \Phi +H(x, \partial_x \Phi)$ vanishes to high order on $\gamma$. The propagation of amplitude can then be determined by $LA_0=0$ to make the $O(\epsilon)$ term vanish. We denote $\Phi(t, x)$ on the curve $\gamma$ by $S$ and the leading amplitude by $A=A_0$. These lead to the standard WKB system
\begin{align}\label{so}
& \partial_t S +H(x, \nabla_x S)=0, \\ \label{ao} & \partial_t A + H_p\cdot \nabla_x A=-
\frac{A}{2}[Tr(H_{xp}) +Tr(\nabla_x^2SH_{pp})],
\end{align}
where $Tr$ is the usual trace map.
We then compute the Taylor series of $ \partial_t \Phi +H(x, \Phi_x)$ about $x(t)$ to first and second order to obtain equations for the phase gradient and Hessian $(u, M)=(\nabla_x \Phi, \nabla_x^2 \Phi)$ as follows
\begin{align}\label{u}
&\partial_t u +H_p\cdot \nabla_x u =-H_x,\\ \label{m} & \partial_t
M+H_p\cdot \nabla_x M +H_{xx} +H_{xp}M+MH_{px} +MH_{pp}M=0.
\end{align}
It is shown in \cite{Ralston82} that the above construction is only possible if $(x(t), p(t))$, where $p(t)=\nabla_x \Phi(t, x(t))$, is a (null) bi-characteristic curve, which is consistent with the characteristic system for the Hamilton-Jacobi equation (\ref{so}).
\begin{align}
\frac{d}{dt}x=H_p, \quad  x(0)=x_0, \\
\frac{d}{dt}p=-H_x, \quad p(0)=p_0,
\end{align}
where $(x, p)=(x, p)(t; x_0, p_0)$. From here on we denote the phase
space variable as $X=(x, p)$ and $X_0=(x_0, p_0)$. Then the
Hamiltonian dynamics can be expressed as
\begin{equation}\label{x}
    \frac{d}{dt}X(t, X_0)=V(X(t, X_0)), \quad X(0, X_0)=X_0.
\end{equation}
The phase velocity $V=(H_p, -H_x)$ is divergence free, i.e.
$div_X(V)=0$.  On this curve $X=X(t, X_0)$, Gaussian beam components
of the first order such as the phase $S(t;X_0)$, the Hessian
$M(t;X_0)$ as well as the amplitude $A(t; X_0)$, are obtained by
solving the following system of ODEs
\begin{align}\label{s}
& \frac{d}{dt} S(t; X_0)=p\cdot H_p-H(x, p), \quad S(0;
X_0)=S_{\rm in}(x_0),\\ \label{m}
 &
 \frac{d}{dt}M(t; X_0) + H_{xx} +H_{xp}M+MH_{px} +MH_{pp}M=0, \quad
 M(0; X_0)=M_{\rm in}(X_0),\\ \label{a}
& \frac{d}{dt} A(t;X_0)  =-\frac{A}{2}\left[ Tr[H_{xp}] +Tr[ M H_{pp}]
\right], \quad A(0; X_0)=A_{\rm in}(x_0).
\end{align}
The essential idea behind the GB method is to choose some complex
Hessian $M_{\rm in}$ initially so that $M$ remains bounded for all time, and its
imaginary part is positive definite. This way the amplitude  $
A(t;X_0)$ is ensured to be also globally bounded from solving
(\ref{a}).

With these components in place, the Gaussian beam phase is constructed
as
$$
\Phi(t,y; X_0)=S(t; X_0)+p(t, X_0)(y-x(t, X_0))
+\frac{1}{2}(y-x(t, X_0))^\top M(t;X_0)(y-x(t, X_0)),
$$
where $p(t, X_0)=\nabla_x \Phi(t, x(t, X_0))$. The leading order of the amplitude is taken as
$$
A(t, y; X_0)=A(t;X_0).
$$
The above construction ensures that the following GB Ansatz is an
approximate solution
$$
\psi_{GB}(t, y; X_0)=A(t;X_0)\exp\left(\frac{i}{\epsilon}\Phi(t,
y;X_0)\right).
$$
The requirement that $Im(M)$ be positive definite ensures that the
asymptotic solution is concentrated on $y=x(t, X_0)$, see e.g. \cite{Ralston82}.

\subsection{Phase space based Gaussian beam ansatz}
If we regard $X_0$ to be the Lagrangian particle marker, then the
map
\begin{equation}\label{xx0}
    X=X(t, X_0)
\end{equation}
serves as a particle trajectory mapping: an initial domain $\Omega
\in R^{2n}$ in phase space evolves in time to
\begin{equation}\label{vp}
    X(t, \Omega)=\{X(t, X_0), \quad X_0\in \Omega \},
\end{equation}
with the vector $V=(H_p, -H_x)$ tangent to the particle trajectory
in phase space. Since the velocity field is divergent-free, the
elementary properties of $X(t, X_0)$ tells that $Vol (X(t, \Omega))=Vol(\Omega)=|\Omega|$
and $$det\left(\frac{\partial X(t, X_0)}{\partial{X_0}}\right)=1.$$
 In other words
the map is volume-preserving and invertible.

The phase space based Gaussian beam ansatz is thus obtained by
changing $X_0$ to $X$ through this particle-trajectory map:
\begin{equation}\label{pgb}
    \psi_{PGB}(t, y, X)=\tilde A(t, X)\exp \left( \frac{i}{\epsilon}
\tilde \Phi(t, y, X) \right),
\end{equation}
where
\begin{equation}\label{tphi}
\tilde \Phi(t, y, X)=\tilde S(t, X)+p\cdot (y-x)
+\frac{1}{2}(y-x)^\top \tilde M(t, X)(y-x).
\end{equation}
To derive the corresponding dynamics for $(\tilde S, \tilde M, \tilde A)$ we  need the following fact.
\begin{lem} \label{lem2.1} ({\bf Operator lifting }) Let the phase
representative of $w(t;X_0)$ be $\tilde w(t, X)$ in the sense that
$w(t;X_0)=\tilde w(t, X(t, X_0)$ for any $t>0$, then
$$
\frac{d}{dt}w(t; X_0) = \mathcal{L} \tilde w(t, X),
$$
where $\mathcal{L}$ is the usual Liouville operator defined by
\begin{equation}\label{liou}
    \mathcal{L} := \partial_t +V\cdot \nabla_X.
\end{equation}
\end{lem}
\begin{proof}
Taking differentiation of
$$
w(t;X_0) \equiv w(t, X(t, X_0)), \quad \forall t>0
$$
in time, we obtain
$$
\frac{d}{dt}w(t;X_0)=\partial_t w+\frac{d}{dt}X(t, X_0)\cdot
\nabla_X w=\partial_t w+V\cdot \nabla_X w.
$$
\end{proof}
Changing the time derivative $\frac{d}{dt}$ to
the Liouville operator $\mathcal{L}$ in the
Lagrangian formulation of equations for $(S, M,
A)$ in (\ref{s})-(\ref{a}), we obtain the  PDEs
for $(\tilde S, \tilde M, \tilde A)$ in phase
space:
\begin{align}\label{sp}
& \mathcal{L}(\tilde S) =p\cdot H_p-H(x, p), \quad \tilde S(0,
X)=S_{\rm in}(x),\\\label{mp}
 &
 \mathcal{L}(\tilde M) + H_{xx} +H_{xp}\tilde M+ \tilde MH_{px} +\tilde MH_{pp} \tilde M=0, \quad
 \tilde M(0, X)=M_{\rm in}(X),\\ \label{ap}
& \mathcal{L}(\tilde A)  =-\frac{\tilde A}{2}\left[ Tr[H_{xp}] +Tr[
\tilde M H_{pp}] \right], \quad \tilde A(0, X)=A_{\rm in}(X).
\end{align}
A method for solving the equation for $\tilde M$
based on Ricatti equations was given by  Leung
and Qian (see \cite{LQ:2008}, formulas (56)-(57))
.  Jin, Wu and Yang have an alternative way of
computing $\tilde M$ based on complex level set
functions (see \cite{JWY:2008}, formulas
(3.5)-(3.16)). We give yet another method of
constructing $\tilde M$, and the higher
derivatives of the phase by level set methods in
\S 7. We point out that though $\psi_{GB}(t,y;
X_0)$ is an asymptotic solution to the wave
equation, $\psi_{PGB}(t, y, X)$ is usually not.
But we shall show that its
integral over the moving domain $X(t, \Omega(0))$ remains
an asymptotic solution of the wave equation. It is this remarkable
feature that allows us to globally recover the original wave field
from only some phase space based measurements!

\section{Recovery of wave fields by superposition}
Since the wave equation we consider is linear, the high frequency wave
field $\psi$ at $(t, y)$ in physical space is expected to be generated by a superposition of neighboring Gaussian beams
\begin{equation}\label{GB+}
\psi^\epsilon(t, y)=Z(n,\epsilon)\int_{\Omega(0)}\psi_{GB}(t,
y;X_0)dX_0,
\end{equation}
where
$$
\Omega(0)=\{X_0, \quad x_0\in {\rm supp}(A_{\rm in}), \quad p_0 \in {\rm
range}(\partial_x S_{\rm in}(x))\}
$$
is an open domain in phase space from which we construct
initial Gaussian beams from the given data.
The normalization parameter $Z(n, \epsilon)$ is
determined by matching the initial data $\psi_0(y)$ so that
$$
\|\psi_0(\cdot) -\psi^\epsilon(0, \cdot)\| \to 0, \quad \epsilon \to
0.
$$

By invoking the volume preserving map $X=X(t, X_0)$ and its inverse $X_0=X_0(t, X)$, we obtain a
phase space based Gaussian beam ansatz
$$
\psi_{PGB}(t, y, X):=\psi_{GB}(t, y; X_0(t, X)).
$$
As remarked earlier, since the map is time dependent, the above phase space based GB ansatz is no longer an
asymptotic solution of the wave equation. We note that their superposition over the moving domain $X(t, \Omega(0))$ remains
a correct asymptotic solution.
\begin{equation}\label{egb}
\psi^\epsilon(t, y)=Z(n,\epsilon)\int_{\Omega(t)}\psi_{PGB}(t,
y,X)dX,
\end{equation}
where
$$
\Omega(t)=X(t, \Omega(0)).
$$
This can be seen directly by using change of variables to go back to the Lagrangian superposition (\ref{GB+}).

In what follows we will construct $\psi_{PGB}$ without reference to  $\psi_{GB}$.  While we could still recover
$\psi_{GB}$ from $\psi_{PGB}$ by coordinate transformations, we can check directly that superpositions of $\psi_{PGB}$
are asymptotic solutions. This only requires the following two lemmas.
\begin{lem}
For any smooth $f(t, X)$ and divergence-free velocity field $V$, one
has
\begin{equation}\label{dt}
    \frac{d}{dt}\int_{X(t, \Omega)}f(t, X)dX=\int_{X(t,
    \Omega)}[\partial_t f +\nabla_X\cdot (fV)]dX.
\end{equation}
\end{lem}
Our estimates are consequences of the following elementary lemma.
\begin{lem}\label{lem3.2}Assume that $Im(\tilde \Phi(t, y, X))\geq c |y-x|^2$, $c>0$, and the Lebeque measure of the initial domain $|\Omega(0)|$
is bounded. Let $B(t, y, X)$ be a smooth function,
satisfying
$$
|B|\leq C|y-x|^k, \quad k>0.
$$
Then we have
$$
\left\|\int_{\Omega(t)} B(t, y, X)e^{i\tilde \Phi(t,y,X)/\epsilon}dX
\right\|_{L^2_y} \lesssim
|\Omega(0)|\ep^{\frac{k}{2}+\frac{n}{4}}.
$$
\end{lem}
\begin{proof}Using Minkowski's integral inequality  we have
\begin{align*}
\left\|\int_{\Omega(t)} B(t, y, X)e^{i\tilde \Phi(t,y,X)/\epsilon}dX
\right\|_{L^2_y} & \leq  \left( \int_{y} \left|
\int_{\Omega(t)}|B|e^{-Im(\tilde \Phi)/\epsilon}dX\right|^2 dy\right)^{1/2} \\
& \leq \int_{\Omega(t)} \left( \int_y |B|^2 e^{-2Im(\tilde
\Phi)/\ep}dy \right)^{1/2}dX \\
& \leq C \int_{\Omega(t)} \left( \int_y |y-x|^{2k} e^{-2c|y-x|^2/\ep}dy \right)^{1/2}dX,
\end{align*}
continuing the estimate with the stretched coordinates $y-x=\ep^{1/2}y'$, and changing from $y$ to $y'$ in the integral
\begin{align*}
& \leq  C \int_{\Omega(t)} \ep^{\frac{k}{2}+\frac{n}{4}}\left( \int_{y'} |y'|^{2k} e^{-2c|y'|^2}dy' \right)^{1/2} dX \\
& =C |\Omega(t)| \ep^{\frac{k}{2}+\frac{n}{4}} \left( \int_{y} |y|^{2k} e^{-2c|y|^2}dy \right)^{1/2},
\end{align*}
which when using $|\Omega(t)|=|\Omega(0)|$ proves the result.
\end{proof}
The normalization parameter needs to be chosen to match the initial data.
For example, if initially $Im(
M_{\rm in})=\beta I$, $\beta>0$, then we need to arrange to match the initial data against $\exp(-\beta|x-y|^2/\ep)$. That accounts for \begin{equation}\label{z}
Z(n, \epsilon)=\left( \int_y e^{-y^\top Im(M_{\rm in})y/(2\epsilon)}
dy\right)^{-1}=\left( \frac{\beta}{2\pi\epsilon}\right)^{n/2}\sim \ep ^{-n/2}
\end{equation}
in dimension $n$. Taking the Schr\"{o}dinger equation as an example, we obtain the
following.
\begin{thm}\label{th3.3} Let $P$ be the linear Schr\"{o}dinger wave operator of the form
$P=-i\epsilon \partial_t +H(y, -i\epsilon \partial_y)$, where $H(y,
p)= \frac{|p|^2}{2} +V_{ext}(y)$, and $\psi^\ep$ is defined in (\ref{egb})
with $Im(\tilde M)$ being positive definite, and $Z(n, \epsilon) \sim \ep^{-n/2}$, then $\psi^\epsilon$ is an asymptotic solution and
satisfies
\begin{equation}\label{pep}
\|P[\psi^\ep](t, \cdot)\|_{L_y^2} \lesssim
|\Omega(0)| \ep^{\frac{3}{2}-\frac{n}{4}}.
\end{equation}
\end{thm}
\begin{proof} We apply the operator $P$  to both sides of (\ref{egb}) to obtain
\begin{align}\notag
Z^{-1} P[\psi^\epsilon ] &= (-i\epsilon \partial_t +H(y, -i\epsilon
\partial_y)) \int_{\Omega(t)}\psi_{PGB}(t, y, X)dX\\ \label{extra}
&= \int_{\Omega(t)} \left[ P[\psi_{PGB}]-i\epsilon \nabla_X\cdot
(V\psi_{PGB})) \right]dX.
\end{align}
By a straightforward calculation it follows
\begin{align*}
&P[\psi_{PGB}]=-i\epsilon \partial_t \tilde A(t, X)e^{i \tilde
\Phi(t, y, X)/\epsilon}+ \tilde A P[e^{i \tilde \Phi(t, y,
X)/\epsilon}]
\\
&=e^{i \tilde \Phi(t, y, X)/\epsilon}\left[ \tilde A [\partial_t
\tilde \Phi+ H(y, \partial_y \tilde \Phi)] - i\epsilon
\left(\partial_t \tilde A +\frac{\tilde A}{2} Tr[\partial_y^2 \tilde
\Phi ]\right)\right].
\end{align*}
The transport term in the integrand gives
$$
-i\epsilon \nabla_X\cdot (V\psi_{PGB})=e^{i \tilde \Phi(t, y,
X)/\epsilon}\left[-i\epsilon V\cdot \nabla_X \tilde A +\tilde A
V\cdot \nabla_X \tilde \Phi \right].
$$
Putting together we have
\begin{equation}\label{pis}
[P-i\ep \nabla_X \cdot V](\psi_{PGB})=e^{i \tilde \Phi(t, y,
X)/\epsilon}\left[ \tilde A [\mathcal{L}[\tilde \Phi]+ H(y,
\partial_y \tilde \Phi)] - i\epsilon \left(\mathcal{L}[\tilde A]
+\frac{\tilde A}{2} Tr[\partial_y^2 \tilde \Phi ]\right)\right].
\end{equation}
Using  $\partial_y^2 \Phi=\tilde M(t, X)$ and that $\tilde A$ solves
(\ref{ap}), we see that $O(\epsilon)$ term vanishes. From
(\ref{tphi}) it follows that
$$
\mathcal{L}[\tilde \Phi]=\mathcal{L}[\tilde S]-|p|^2
-\partial_x V_{ext}\cdot(y-x)+\frac{1}{2}(y-x)^\top \mathcal{L}[\tilde
M](y-x)-p^\top \tilde M(y-x)
$$
and
$$ H(y, \partial_y \tilde \Phi)=V_{ext}(y)+\frac{1}{2}|p+\tilde
M(y-x)|^2.$$ This when using equation (\ref{sp}) for $\tilde S$
gives
\begin{equation}\label{yx3}
\mathcal{L}[\tilde \Phi]+ H(y,
\partial_y \tilde \Phi)=V_{ext}(y)-V_{ext}(x)- \partial_x V_{ext}(y-x)+\frac{1}{2}(y-x)^\top
\partial^2_xV_{ext}(y-x)=O(|y-x|^3).
\end{equation}
Consequently, using the above lemma with $k=3$,
\begin{align*}
\|P[\psi^\ep](t, \cdot)\|_{L^2_y}\lesssim Z(n, \epsilon)
|\Omega(0)|\ep^{\frac{3}{2}+\frac{n}{4}},
\end{align*}
which with (\ref{z}) leads to the desired estimate.
\end{proof}
\begin{rem}
The extra term in the integral in (\ref{extra}) gives an alternate way of seeing that the phase space super-position is an accurate solution of the PDE. Of course, it integrates to zero when the support of the beam superposition does not touch the boundary of the integration domain, but it makes it possible to verify the accuracy without going back to the Lagrangian super-position.
\end{rem}

We now obtain the following  estimate.
\begin{thm}
Given $T>0$, and let $\psi$ be the solution of the Schr\"{o}dinger
equation subject to the initial data $\psi_0$, and $\psi^\epsilon$
be the approximation defined in (\ref{egb}) with
$Im(M_{\rm in})$ being positive definite, and $|\Omega(0)|<\infty$. Then
there exists $\epsilon_0>0$, a normalization parameter $Z(n,
\epsilon) \sim \ep ^{-n/2}$, and a constant $C$ such that for all $\epsilon \in (0,
\epsilon_0)$
$$
\|(\psi^\epsilon -\psi)(t, \cdot)\|_{L^2}\leq \|\psi^\epsilon(0,
\cdot)-\psi_{\rm in}(\cdot) \|_{L^2} +C|\Omega(0)|\ep^{\frac{1}{2}-\frac{n}{4}}
$$
for $t\in [0, T]$.
\end{thm}
\begin{proof}Let $e:=\psi^\epsilon-\psi$, then from $P[\psi]=0$ it follows
$$
P[e]=P[\psi^\ep]-P[\psi]=P[\psi^\ep].
$$
A calculation of $\int_{R^n} [e \overline{P[e]}-\bar e P[e]]dy$
leads to
$$
\epsilon \frac{d}{dt}\int_y |e|^2dy =\int_y Im(e\overline{P[\psi^\ep]})dy.
$$
Integration over $[0, t]$ gives
\begin{equation}\label{well}
\|e(t, \cdot)\|_{L^2_y}\leq \|e(0,
\cdot)\|_{L_y^2}+\frac{1}{\epsilon}\int_0^t \|P[\psi^\ep](\tau,
\cdot)\|_{L^2}d\tau, \quad t\in [0, T].
\end{equation}
This when combined with the estimate for $P[\psi^\ep]$ in (\ref{pep})  gives the result
as desired.
\end{proof}
\begin{rem}
The approximation error comes from two sources: initial error and
the evolution error. To improve accuracy one has to enhance the
accuracy for both. The evolution accuracy can be improved by
obtaining more phase space based measurements such as higher order
derivatives of phase and amplitude, which will be sketched in \S 5.
\end{rem}
\begin{rem}
In phase space the tracking of the beam propagation is lost in the
Gaussian beam ansatz, but has been recorded through the moving
domain $\Omega(t)$, which can be traced back to $\Omega(0)$.
\end{rem}
\section{Control of the initial error}
Let $K(x, \tau)=\frac{1}{(4\pi \tau)^{n/2}}e^{-\frac{|x|^2}{4\tau}}$
be the usual heat kernel, satisfying ${\rm limit}_{\tau \downarrow
0}K(x, \tau)=\delta (x)$ as distributions on $\mathbb{R}^n$. Then
$$
\int_x K(x-y, \tau)dx=1, \quad \forall \tau>0, \; y\in \mathbb{R}^n.
$$
For highly oscillatory initial data we have
$$
\psi_{\rm in}(y)=A_{\rm in}(y)e^{iS_{\rm in}(y)/\epsilon}=\int_x
A_{\rm in}(y)e^{iS_{\rm in}(y)/\epsilon}K\left(x-y, \frac{\epsilon}{2}\right)dx.
$$
Both the phase and amplitude in the integrand can be approximated by their Taylor expansion
when $|x-y|$ is small, say  $|x-y|<\epsilon^{1/3}$, and the integral will then be
$O(\exp (-c/\epsilon^{1/3}))$ with some $c<\frac{1}{2}$ outside this neighborhood.
Let $T_j^x[f](y)$ denote the $j^{th}$ order Taylor polynomial of $f$ about $x$ at the point $y$. Then
\begin{equation}\label{ip}
    \psi_{\rm in}(y)\sim \int_x A_{\rm in}(x)e^{\frac{i}{\epsilon}\left[T_{2}^x[S_{\rm in}](y)\right]} K\left(x-y, \frac{\epsilon}{2}\right)dx,
\end{equation}
which tends to $\psi_{\rm in}$ as $\epsilon \to 0$.

Indeed, the approximate accuracy is ensured by the following result by Tanushev \cite{Tan08}.
\begin{lem} \label{tan}
Let $S_{\rm in} \in C^\infty (R^n)$ be a real-valued function, and $A_{\rm in} \in
C_0^\infty (R^n) $, and $\rho \in C_0^\infty (R^n)$ be such that
$\rho\geq 0$, $\rho \equiv 1$ in a ball of radius $\delta  > 0$
about the origin.
Define,
\begin{align*}
\psi_{\rm in}(y) & = A_{\rm in}(y)e^{iS_{\rm in}(y)/\epsilon},\\
v(y; x) &=\rho(y-x)T_j^x[A_{\rm in}](y)e^{\frac{i}{\epsilon}\left[T_{j+2}^x[S_{\rm in}](y)\right]}
K\left(x-y, \frac{\epsilon}{2}\right).
\end{align*}
Then
$$
\left\|\psi_{\rm in}(\cdot)-\int_{supp(A_{\rm in})}v(\cdot; x)dx \right\|_{L^2} \lesssim
\ep^{\frac{j+1}{2}}.
$$
\end{lem}
\begin{rem}We note that for the case $j=0$ the cutoff function $\rho$ is unnecessary, and it can be shown the initial approximation error remains of order $\ep^{1/2}$. However, a cutoff function is
certainly important when one is building beams of higher accuracy because the higher order terms in the Taylor expansion of the phase can change the sign of its imaginary part when one does not
stay close to the central ray, see \S 5.
\end{rem}

We take the above approximation in (\ref{ip}) as initial data for $\psi^\ep(0, y)$
and rewrite it as follows
\begin{equation}\label{delta}
 \psi^\epsilon(0, y)=Z(n, \epsilon) \int_{\Omega(0)}\psi_{PGB}(0,y,
X)\delta(p-\nabla_x S_{\rm in}(x))dX, \quad Z(n, \ep)=\frac{1}{(2\pi \ep)^{n/2}}
\end{equation}
with $\tilde A_{\rm in}=A_{\rm in}(x), \; \tilde S_{\rm in}(X)=S_{\rm in}(x)$ and $\tilde
M_{\rm in}(X)=\partial_x^2 S_{\rm in}(x)+i I$. The above lemma ensures that
\begin{equation}\label{e0}
\left\|\psi_{\rm in}(\cdot)-\psi^\ep(0, \cdot) \right\|_{L^2} \lesssim \ep^{1/2}.
\end{equation}
In order to track the deformation of the surface $p-\nabla_x S_{\rm in}(x)=0$ as time evolves, we introduce a function $w=w(t, X)$ such that
\begin{equation}\label{w}
\mathcal{L}[w]=0, \quad w(0, X)=p-\nabla_x S_{\rm in}(x).
\end{equation}
For smooth Hamiltonian $H(x, p)$, $w$ remains smooth once it is initially so.  A modified  approximation is defined as \begin{equation}\label{psw}
\psi^\ep(t, y):=Z(n, \ep)\int_{\Omega(t)}\psi_{PGB}(t, y, X)\delta(w(t, X))dX,
\end{equation}
which has taken care of the Dirac delta function in (\ref{delta}).  We then have the following theorem.
\begin{thm}
If assumptions of Theorem \ref{th3.3} are met, then $\psi^\ep$ defined in (\ref{psw}) is also an asymptotic solution, satisfying
\begin{equation}\label{ep2}
   \|P[\psi^\ep](t, \cdot)]\|_{L^2} \lesssim |{\rm supp}(A_{in})|\epsilon^{\frac{3}{2}-\frac{n}{4}}.
\end{equation}
\end{thm}
\begin{proof}Using the volume-preserving map of $X=X(t, X_0)$ and $w(t, X(t, X_0)=w(0, X_0)$, we have
\begin{align}\notag
\psi^\epsilon(t, y)&  =Z(n,\epsilon)\int_{\Omega(0)}\psi_{PGB}(t,
y, X(t, X_0))\delta(w(t, X(t, X_0)))dX_0 \\ \notag
& =Z(n,\epsilon)\int_{\Omega(0)}\psi_{GB}(t,
y; X_0)\delta(w(0, X_0))dX_0\\ \notag
&= Z(n,\epsilon)\int_{\Omega(0)}\psi_{GB}(t,
y; X_0)\delta(p_0-\nabla_x S_{\rm in}(x_0))dX_0\\ \label{gbs}
& =Z(n,\epsilon)\int_{{\rm supp}(A_{\rm in})} A(t; x_0) G(t, y; x_0)e^{i\Phi(t, y;x_0)/\ep}dx_0.
\end{align}
Here for simplicity we use only $x_0$  in the integrand instead of $X_0=(x_0, \nabla_x S_{in}(x_0))$.
We now repeat the similar estimate to that in the proof of Lemma \ref{lem3.2} with $k=3$ to obtain
\begin{equation}\label{xx}
\left\| \int_{{\rm supp}(A_{\rm in})} A(t; x_0) G(t, y; x_0)e^{i\Phi(t, y;x_0)/\ep}dx_0 \right\|_{L^2_y}
\lesssim |{\rm supp}(A_{\rm in})|\ep^{\frac{3}{2}+\frac{n}{4}}.
\end{equation}
Hence
\begin{align*}
\|P[\psi^\ep](t, \cdot)\|_{L^2}\lesssim Z(n, \epsilon)
||{\rm supp}(A_{\rm in})|\ep^{\frac{3}{2}+\frac{n}{4}}
\end{align*}
which with (\ref{z}) leads to the desired estimate.
\end{proof}
Plugging estimates (\ref{e0}) and (\ref{ep2}) into (\ref{well})
we arrive at our main result.
\begin{thm}\label{thm4.3}
Given $T>0$, and let $\psi$ be the solution of the Schr\"{o}dinger
equation subject to the initial data $\psi_{\rm in}=A_{\rm in}e^{iS_{\rm in}(x)/\ep}$, and $\psi^\epsilon$
be the first order approximation defined in (\ref{psw}) with initial data
satisfying $\tilde S_{\rm in}(X)=S_{\rm in}(x)$,  $\tilde M_{\rm in}(X) =
\partial_x^2 S_{\rm in}(x)+i I$, and $\tilde A_{\rm in}(X)=A_{\rm in}(x)$ with $|supp(A_{\rm in})|<\infty$.
Then there exists $\epsilon_0>0$, a normalization parameter $Z(n,
\epsilon)$, and a constant $C$ such that for all $\epsilon \in (0,
\epsilon_0)$
$$
\|(\psi^\epsilon -\psi)(t, \cdot)\|_{L^2} \lesssim |{\rm supp}(A_{\rm in})|\ep^{\frac{1}{2}-\frac{n}{4}}
$$
for $t\in [0, T]$.
\end{thm}
\begin{rem} The exponent $1/2$ reflects the accuracy of the Gaussian beam in solving the PDE. It will increase when one uses more accurate beams. 
The exponent $-\frac{n}{4}$ indicates the blow-up rate for the worst possible case due to caustics.
Of course, if nature of the caustic were a priori known, it would be possible to obtain a better convergence rate by taking the caustic structure into account.
\end{rem}

\section{A closer look at caustics}
\subsection{Schur's lemma}
Instead of using the Minkowski inequality we shall use Schur's lemma to see how caustic structure may be used to obtain a better estimate. Recall Schur's Lemma: If $[Tf](y) =\int
K(x, y)f(x)dx$ and
$$
{\rm sup}_x
\int_y |K(x, y)|dy = C_1, \; {\rm sup}_y
\int_x |K(x, y)|dx = C_2,
$$
then
$$
\|Tf\|_{L^2}\leq \sqrt{C_1C_2}\|f\|_{L^2}.
$$
\begin{proof}
We have by Schwartz
\begin{align*}
|[Tf](y)|^2 \leq
\left(
\int|K(x,y)|f(x)dx\right)^2 & \leq
\int |K(x,y)|dx\int
|K(x,y)||f(x)|^2dx\\
 & \leq  C_2 \int
|K(x,y)||f(x)|^2dx.
\end{align*}
So integrating both sides in $y$ and taking the square root gives the result.
\end{proof}
We now apply Schur's lemma to left hand side of (\ref{xx}).  So for simplicity
$$
[Tf](y) = \int_{{\rm supp}(A_{\rm in})} A(t; x_0) G(t, y; x_0)e^{i\Phi(t, y;x_0)/\ep}dx_0,
$$
where the imaginary part of $\Phi(t, y;x_0)$ is bounded below by $cI$ and for convenience we will assume that $|G|\leq |y-x(t, x_0)|^k$. Then one can apply Schur's lemma with
$$
C_1={\rm sup}_{x_0}\int_y |y-x(t, x_0)|^k e^{-(c/\ep)|y-x(t, x_0)|^2}dy=\ep^{\frac{k}{2}+\frac{n}{2}}\int_z |z|^ke^{-c|z|^2}dz, \quad {\rm and}
$$
$$
C_2(t, \ep)={\rm sup}_{y}\int_{x_0}|y-x(t, x_0)|^k e^{-(c/\ep)|y-x(t, x_0)|^2}dx_0.
$$
In general one does not know what $C_2(t, \ep)$ will be. As long as $A$ has compact support $C_2$ will be at least
bounded by $c\ep^{k/2}$.  Thus the error in $L^2$ norm will be bounded by $c\ep^{k/2 +n/4}$,  as shown in (\ref{xx}) with the Minkowski inequality.
Note that the worst case for the wave equation, $\partial_t^2 \psi -\Delta \psi=0$, occurs when $x(t, x_0) = x_0(1- t/|x_0|)$. In that case $C_2=c \ep^{k/2 +(1+n)/4}$, yielding a better rate.

\subsection{An example with remarkable accuracy}
Here below we illustrate that a better convergence rate can also be obtained for the Schr\"{o}dinger equation with quadratic potential and quadratic phase. We consider the solution of
$$
i\ep \partial_t \psi  = -\frac{\ep^2}{2} \Delta \psi, \quad  x\in \mathbb{R}^n,
$$
with the initial data $\psi(0, x) = \exp(-i|x|^2/(2\ep))$, then
$$
\psi(t, x) = (1-t)^{-n/2} \exp \left(-\frac{i|x|^2}{2\ep (1-t)}\right).
$$
This solution becomes a multiple of the $\delta$-function at $t =1$. This suggested solving
the free Schr\"{o}dinger equation with initial data $\psi(x,0) = g(x)\exp\left(-i|x|^2/(2\ep)\right)$ where $g\in C^{\infty}(\mathbb{R}^n)$ and evaluating the
solution at $t =1$.  Using Fourier transform one may express the solution  as
$$
\psi(t, x)=\frac{1}{(2\pi i \ep t)^{n/2}}\int_y \psi_0(y)\exp \left(\frac{i}{2\ep t}|x-y|^2 \right)dy.
$$
Evaluating  at $t=1$, we have
$$
\psi(1, x)=\frac{1}{(2\pi i \ep )^{n/2}}\int_y g(y)e^{-ix\cdot y/\ep}dy e^{i|x|^2/(2\ep)}=c\ep^{-n/2} \hat g( x/\ep)e^{i|x|^2/(2\ep)},
$$
where $c=e^{-\pi n i/4}$ and $\hat g$ is the Fourier transform of $g$, defined by
$$
\hat g(\xi)=\frac{1}{(2\pi)^{n/2}}\int g(x) e^{-ix\cdot \xi}dx.
$$
It is easy to verify that
$$
\|\psi(1, \cdot)\|_{L^2} =\|\psi(0, \cdot)\|_{L^2}.
$$
As $\ep \to 0$, $\psi(1, x)$ diverges (pointwise) like $\ep^{-n/2}$  near $x = 0$, but goes rapidly to
zero away from $x = 0$.

We now build a superposition of Gaussian beams approximation for the solution of the same problem. For the Gaussian beam superposition, the phase is obtained as
$$
\Phi(t, x; y)=(t-1)\frac{|y|^2}{2}- x\cdot y+|y|^2(1-t) +\frac{\beta i-1}{1+(\beta i-1)t}\frac{|x-y(1-t)|^2}{2}.
$$
Note that we have chosen $\beta >$ 0 for the initial beam width. For the amplitude we get
$$
(1 + (\beta i-1)t)^{n/2}A(t, y(1-t)) = A(0, y) = g(y).
$$
So, setting $x(t; y) = (1-t)y$, we obtain
$$
A(t, x(t; y)) = (1 +(\beta i -1)t)^{-n/2}g(y).
$$
If we do the superposition with the normalization, we end up with
$$
\psi^\ep(t, x)= \left(\frac{\beta }{2\pi \ep}\right)^{n/2}
\int_y [1+(\beta i -1)t]^{-n/2}g(y)e^{i\Phi(t, x;y)/\ep}dy.
$$
If we evaluate that at $t =1$, it becomes
\begin{align*}
\psi^\ep(1, x) &= \left(\frac{\beta }{2\pi \ep}\right)^{n/2}
\int_y [\beta i ]^{-n/2}g(y)e^{\frac{i}{\ep}[- x\cdot y + \frac{\beta i-1}{\beta i}\frac{|x|^2}{2}]}dy\\
&=c\ep^{-n/2}\hat g(x/\ep) e^{(\beta i-1)|x|^2/(2\beta\ep)}\\
&=\psi(1, x)e^{-|x|^2/(2\beta \ep)}.
\end{align*}
This shows that at the caustic $x=0$, both become the same. We can see the error $\psi(1, x)\left(1-e^{-|x|^2/(2\beta \ep)}\right)$ when measured in $L^2$-norm:
\begin{align*}
\|\psi^\ep(1, \cdot)-\psi(1, \cdot)\|_{L^2}^2 &=\ep^{-n}\int |\hat g\left(\frac{x}{\ep}\right)|^2
\left(1-e^{-|x|^2/(2\beta \ep)}\right)^2dx\\
&=\int |\hat g\left(z\right)|^2
\left(1-e^{-\ep|z|^2/(2\beta )}\right)^2dz.
\end{align*}
That implies that \\
(a) For any $g\in L^2$ the Gaussian beam approximation converges to the true solution
(at $t =1$), but there is no uniform estimate on the difference in terms of the
$L^2$-norm of $g$ (an example of strong but not uniform convergence). \\
(b) If
$\int |\hat g(z)|^2(1 + |z|^2)^2dz < \infty$, i.e. if $g\in H^2$, then the norm of the difference is
$O(\ep)$. 

Actually in the current example, it can be verified that the evolution error is zero, so the initial error should propagate
in time. If we look at the initial error of the Gaussian beam approximation, we have
\begin{align*}
\|\psi^\ep(0, \cdot)-\psi(0, \cdot)\|_{L^2}^2 &=\int_x \left | \left(\frac{\beta }{2\pi \ep}\right)^{n/2} \int_y g(y)
e^{-\frac{\beta}{2\ep}|x-y|^2}dy -g(x) \right|^2 dx\\
&=\left(\frac{\beta }{2\pi \ep}\right)^{n}  \int_x \left| \int_y (g(y)-g(x))
e^{-\frac{\beta}{2\ep}|x-y|^2}dy\right|^2 dx.
\end{align*}
Set $K(x)=\left(\frac{\beta }{2\pi \ep}\right)^{n/2}e^{-\frac{\beta}{2\ep}|x|^2}$, a direct integration shows that $\hat K=(2\pi)^{-n/2} e^{-\ep|\xi|^2/(2\beta )}$. If we apply Parseval's theorem, we obtain
\begin{align}\notag
\|\psi^\ep(0, \cdot)-\psi(0, \cdot)\|_{L^2}^2 &=\| \hat g -(2\pi)^{n/2} \hat g\cdot \hat K\|_{L^2}^2\\ \label{ii+}
& =\int |\hat g\left(z\right)|^2
\left(1-e^{-\ep|z|^2/(2\beta )}\right)^2dz,
\end{align}
which is the same error as that evaluated at $t=1$. We point out that for the initial phase of general form, the initial error is still $O(\ep^{1/2})$ unless an higher order expansion of the phase is used.

There are two conclusions that one can draw from the preceding.
\begin{thm}Under assumptions of Theorem \ref{thm4.3} and assume that the potential is a quadratic function. Then for $t\in[0, T]$ and $\ep \in (0, \ep_0)$ we have
\begin{itemize}
\item If $S_{\rm in}$ is a quadratic function and $A_{\rm in} \in H^2$
$$
\|(\psi^\epsilon -\psi)(t, \cdot)\|_{L^2} \lesssim   \ep.
$$
\item If $S_{\rm in}\in C^\infty$ and  $A_{\rm in} \in C_0^\infty$
 $$
\|(\psi^\epsilon -\psi)(t, \cdot)\|_{L^2} \lesssim   \ep^{1/2}.
$$
\end{itemize}
\end{thm}
\begin{proof}It follows from (\ref{yx3}) that for quadratic potentials
$$
P[\psi^\ep]=0.
$$
Then the total error is governed by the initial error only. For quadratic potentials and $A_{\rm in}\in H^2(\mathbb{R}^n)$
we obtain the $O(\ep)$ error as shown in (\ref{ii+}). For the general phase function the claim follows from Lemma \ref{tan} with $j=0$.
\end{proof}

\section{Higher order Approximations}
The accuracy of the phase space based Gaussian beam superposition also
depends on accuracy of the individual Gaussian beam ansatz. If we refer the above construction as a first order GB solution, then
the $k^{th}$ order GB solution should involve
terms up to $(k+1)^{th}$ order for the phase, and $(k-1-2l)^{th}$ order for the $l^{th}$ amplitude $A_l$
for $l=0, \cdots, \left[ \frac{k-1}{2}\right]$. The equations for these phase and amplitude Taylor coefficients
can be derived by letting the leading order ones to hold on $\gamma$ along with several of their derivatives (see the Appendix ).

Let $X=X(t; X_0)$, with $x=x(t;X_0)$, denote the characteristic path at time $t>0$, which originates from $X_0$. Following \cite{Tan08} we define the $k^{th}$ order Gaussian beams as follows
$$
\psi_{kGB}(t, y; X_0)=\rho(y-x)\left[ \sum_{l=0}^{\lfloor \frac{k-1}{2}\rfloor }
\ep^l T_{k-1-2l}^{x}[A_l](y) \right] \exp\left(\frac{i}{\epsilon}T_{k+1}^{x}[\Phi](y)\right),
$$
where $T_k^x[f](y)$ is the $k^{th}$ order Taylor polynomial of $f$ about $x$ evaluated at $y$, and $\rho$ is a cut-off function such that on its support the Taylor expansion of $\Phi$ still has a positive imaginary part.

By invoking the volume preserving map $X=X(t, X_0)$ and its inverse map denoted by $X_0=X_0(t, X)$,  we obtain a
phase space based $k^{th}$ order Gaussian beam Ansatz
$$
\psi_{kPGB}(t, y, X):=\psi_{kGB}(t, y; X_0(t, X)).
$$

Beyond the first order GB components, all Taylor coefficients $\partial_y^\alpha \Phi$  for $|\alpha| \geq 3$ in  phase space are replaced by $m_\alpha(t, X)$, satisfying a linear equation (\ref{liu+}) in phase space; and Taylor coefficients $\partial_y^\alpha A_l$  for $|\alpha| \geq 1$ in the amplitude are replaced by $\tilde A_{l, \alpha}$, which can be obtained recursively from solving transport equations in phase space, see the appendix for details.

Proceeding as previously, we form the superpositions.
\begin{equation}\label{kegb}
\psi_k^\epsilon(t, y)=Z(n,\epsilon)\int_{\Omega(t)}\psi_{kPGB}(t,
y,X)\delta(w(t, X))dX,
\end{equation}
where $
\Omega(t)=X(t, \Omega(0)),
$
and $w(t, X)$ is the solution of the Liouville equation subject to $w(0, X)=p-\nabla_x S_{\rm in}(x)$.

This gives a $k^{th}$ order asymptotic solution of the wave equation. More precisely, we have the following theorem.
\begin{thm}\label{3.3} Let $P$ be the linear Schr\"{o}dinger wave operator of the form
$P=-i\epsilon \partial_t +H(y, -i\epsilon \partial_y)$, where $H(y,
p)= \frac{|p|^2}{2} +V_{ext}(y)$, and $\psi^\ep$ is defined in (\ref{kegb})
with $Im(M_{\rm in})=I$ and $Z(n, \epsilon)=(2\pi\ep)^{-n/2}$, then $\psi_k^\epsilon$ is an asymptotic solution and
satisfies
\begin{equation}\label{pk}
   \|P[\psi_k^\ep](t, \cdot)\|_{L_y^2} \lesssim
|{\rm supp}(A_{\rm in})|\ep^{\frac{k}{2}+1-\frac{n}{4}}.
\end{equation}
\end{thm}
\begin{proof}Using the volume-preserving map of $X=X(t, X_0)$ and $w(t, X(t, X_0))=w(0, X_0)$, we have
\begin{align*}
\psi_k^\epsilon(t, y)&  =Z(n,\epsilon)\int_{\Omega(0)}\psi_{kPGB}(t,
y, X(t, X_0))\delta(w(t, X(t, X_0)))dX_0\\
& =Z(n,\epsilon)\int_{\Omega(0)}\psi_{kGB}(t,
y; X_0)\delta(w(0, X_0))dX_0\\
&= Z(n,\epsilon)\int_{\Omega(0)}\psi_{kGB}(t,
y; X_0)\delta(p_0-\nabla_x S_{\rm in}(x_0))dX_0\\
& =Z(n,\epsilon)\int_{{\rm supp} \{A_{\rm in}\}}\psi_{kGB}(t,
y; x_0)dx_0.
\end{align*}
According to the GB construction sketched in the appendix, $\psi_{kGB}(t,
y; x_0)$ are asymptotic solutions for each $x_0$, so will be their superpositions $\psi_k^\epsilon(t, y)$.
It remains to verify (\ref{pk}).  First we see that
$$
P[\psi_k^\epsilon(t, y)]=Z(n,\epsilon)\int_{{\rm supp} \{A_{\rm in}\}}P[\psi_{kGB}(t,
y; x_0)]dx_0.
$$
Using (\ref{pap}) in the appendix with $A$ replaced by $\rho(y-x)\left[ \sum_{l=0}^{\lfloor \frac{k-1}{2}\rfloor }
\ep^l T_{k-1-2l}^{x}[A_l](y) \right]$ and $\Phi$ by
$T_{k+1}^x[\Phi](y)$, we have
$$
c_0(t, y)=[\partial_t T_{k+1}^x[\Phi](y)+H(y, \nabla_y T_{k+1}^x[\Phi](y))]\rho(y-x)T_{k-1}^x[A_0](y).
$$
Using $T_{k+1}^x[\Phi](y)=\Phi(y)+R_{k+1}^x[\Phi](y)$, here $R^x_{k+1}$ denotes the remainder of the Taylor expansion,  and $G(t, y)=\partial_t \Phi +H(y, \nabla_y \Phi)=O(|y-x|^{k+2})$ we can see that
$$
|c_0(t, y)|\leq C|y-x|^{k+2}.
$$
Also using the construction for $A_l$ and their derivatives in the appendix, we are able to show
$$
|c_{l}(t, y)|\leq C|y-x|^{k+2-2l},
$$
where we have used the fact that differentiation of $\rho$ vanishes in a neighborhood of $x$.   The use of the cut-off function ensures that we can always choose a small neighborhood of $x(t; x_0)$ so that
$$
Im(T_{k+1}^x[\Phi](y))\geq c |y-x|^2.
$$
Consequently, using Minkowski's integral inequality,
\begin{align*}
 Z^{-1} \|P[\psi_k^\epsilon(t, \cdot)]\|_{L^2} & \leq  \left( \int_y \left|
\int_{ {\rm supp}(A_{\rm in})} e^{-Im(T_{k+1}^x[\Phi](y))/\epsilon}\left|c_0+c_1\ep + \cdots \right|dx_0 \right|^2 dy\right)^{1/2} \\
& \leq   \int_{{\rm supp}(A_{\rm in})} \left( \int_y
 e^{-2c|y-x(t, x_0)|^2/\epsilon}\left|c_0+c_1\ep + \cdots  \right|^2 dy \right)^{1/2}dx_0 \\
& \leq  C \int_{{\rm supp}(A_{\rm in})}  \left( \int_y
 e^{-2c|y-x(t, x_0)|^2/\epsilon} \sum_{l=0}^{\lfloor \frac{k-1}{2}\rfloor } |y-x(t, x_0)|^{2(k+2-2l)}\ep^{2l}  dy \right)^{1/2}dx_0.
\end{align*}
If we introduce the stretched coordinates $y-x(t;x_0)=\ep^{1/2}y'$, and changing from $y$ to $y'$ in the integral, we see that the new integrand is bounded by
$$
\ep^{k+2+\frac{n}{2}}|y'|^{2(k+2-2l)}\exp\left(-2c|y'|^2\right).
$$
Thus  $\|P[\psi_k^\epsilon(t, \cdot)]\|_{L^2}$ is bounded by $Z(n, \ep)|{\rm supp}(A_{\rm in})| \ep^{\frac{k}{2}+1+\frac{n}{4}}$. The desired estimate then follows.
\end{proof}
In order to obtain an estimate of $\|(\psi_k^\ep -\psi)(t, \cdot)\|$ for any $t\leq T$, all that remains to verify is that  the superposition (\ref{kegb}) accurately approximates the initial data. For $t=0$, the approximation is as follows
\begin{align*}
\psi_k^\ep(0, y) &=Z(n,\ep)\int_{\Omega(0)}\psi_{kPGB}(0,
y,X)\delta(w(0, X))dX \\
& =Z(n,\ep)\int_{\Omega(0)}\psi_{kGB}(0,
y,X_0)\delta(p_0-\nabla_x S_{\rm in}(x_0)))dX_0\\
& =Z(n,\ep)\int_{{\rm supp}A_{\rm in}}\psi_{kGB}(0, y, x_0, \nabla_xS_{\rm in}(x_0)) dx_0,
\end{align*}
where
$$
\psi_{kGB}(0, y, x, \nabla_xS_{\rm in}(x))=
\rho(y-x)\left[
T_{k-1}^{x}[A_{\rm in}](y) \right] \exp\left(\frac{i}{\epsilon}T_{k+1}^{x}[S_{\rm in}](y)\right)e^{-i|y-x|^2/(2\ep)},
$$
 where we have taken $A_0=A_{\rm in}$ and $A_l=0$ for $l\geq 1$,  $\partial_x^\alpha \Phi(0, x)=\partial_x^\alpha S_{\rm in}(x) (\alpha \not=2)$, and
$\partial_x^2 \Phi(0, x)=\partial_x^2 S_{\rm in}(x)+i I.$ From Lemma \ref{tan}  we have that
$$
\left\|\psi_{\rm in}-\psi^\ep(0, \cdot) \right\|_{L^2} \lesssim
\ep^{\frac{k}{2}}.
$$
Thus our main result for $k^{th}$ order phase space GB superposition is as follows.
\begin{thm}
Given $T>0$, and let $\psi$ be the solution of the Schr\"{o}dinger
equation subject to the initial data $\psi_{\rm in}=A_{\rm in}e^{iS_{\rm in}(x)/\ep}$, and $\psi^\epsilon$
be the $k^{th}$ order approximation defined in (\ref{kegb}) with initial data chosen as described above with $|supp(A_{\rm in})|<\infty$.
Then there exists $\epsilon_0>0$, a normalization parameter $Z(n,
\epsilon)\sim \ep^{-n/2}$, and a constant $C$ such that for all $\epsilon \in (0,
\epsilon_0)$
$$
\|(\psi^\epsilon -\psi)(t, \cdot)\|_{L^2} \lesssim |{\rm supp}(A_{\rm in})|\ep^{\frac{k}{2}-\frac{n}{4}}
$$
for $t\in [0, T]$.
\end{thm}

\section{Computing Taylor coefficients of the phase via level set functions}
We now turn to construction of the phase space ingredients required
for the approximation. In order to identify a bi-characteristic
curve in phase space, we  introduce a vector-valued level
set function $\phi \in \mathbb{R}^{2n}$ so that the interaction of zeros of each component
uniquely defines the target curve.  In other words, we
assume that
$$\Gamma=\{(t, X), \; \phi(t, X)= \phi(0, X_0)\} $$
contains the bi-characteristic curve starting from $X_0=(x_0, p_0)$
for any $t>0$, then  $\phi$ must satisfy
$$
\phi(t, X(t, X_0))\equiv \phi(0, X_0).
$$
This is equivalent to the following Liouville equation
\begin{equation}\label{li}
    \mathcal{L}[ \phi (t, X)]=0,
\end{equation}
where $\mathcal{L}:=\partial_t+V\cdot \nabla_X $ is the Liouville
operator. The initial data can be simply taken as
\begin{equation}\label{phi12}
\phi(0, X)=X-X_0.
\end{equation}
Then the curve $\Gamma$ is globally determined by the zero set of a
vector level set function $\phi=(\phi_1, \phi_2)^\top$.

For the construction,  $\tilde S$ can be solved from (\ref{sp}), we
are then left to determine $\tilde M$, followed by solving
(\ref{ap}) to obtain  $\tilde A$. Note that equation (\ref{mp}) is
nonlinear in $\tilde M$, the solution might not exist for all $t>0$.
The heart of the GB method is to choose complex initial data so
that a global
solution is guaranteed and satisfies two requirements \cite{Ralston82}:\\
i) $\tilde M=\tilde M^T,$\\
ii)$Im(\tilde M)$ must be positive definite for all $t>0$.

\subsection{Evaluation of the Hessian}
We now show this can be done via the obtained level set functions
$\phi \in \mathbb{R}^{2n}$.

\begin{thm}\label{thm2.3}
Let $\phi=(\phi_1, \phi_2)^\top$ with $\phi_i \in \mathbb{R}^n$ be the global solution of (\ref{li}) with
the initial condition (\ref{phi12}).  We have \\
a) $\mathcal{L}(k_1\phi_1 +k_2\phi_2)=0$ for any $k_1, k_2\in \mathbb{C}$.\\
b) Set $g:=k_1\phi_1+k_2 \phi_2$. If  $Im(\bar k_1 k_2)\not=0$, then $g_p$ is invertible for all $t>0$.\\
c) If $M_{\rm in}=-g_x(g_p)^{-1}|_{t=0}$, then $\tilde M=-g_x
(g_p)^{-1} $ for all $t>0$,\\
d)$M$ satisfies (\ref{mp}) and i). If $Im(k_1/k_2)<0$, then $\tilde
M$ satisfies ii) too.
\end{thm}
\begin{proof}
a) This follows by noting that the Liouville operator is linear and all its coefficients are real. \\
b) By taking the gradients $\nabla_{x}$ and $\nabla_{p}$ of the Liouville equation $\mathcal{L}(g)=0$, respectively, we obtain the following equations
\begin{align}\label{p1}
\mathcal{L}(g_x)& =H_{xx} g_p- H_{xp} g_x,\\ \label{p2}
\mathcal{L}(g_p)& = H_{px} g_p- H_{pp} g_x.
\end{align}
The equation is understood to be satisfied by each matrix. Let
$B=\overline{g_p}^T g_x-\overline{g_x}^T g_p$ be a
complex matrix, and $I$ an identity matrix, a direct verification
shows that
\begin{align*}
\mathcal{L}(B) & =\mathcal{L}(\overline{g_p}^T g_x)- \mathcal{L}(\overline{g_x}^Tg_p)\\
               & =\overline{\mathcal{L}(g_p)}^T g_x + \overline{g_p}^T\mathcal{L}(g_x)
               - \overline{\mathcal{L}(g_x)}^T g_p + \overline{g_x}^T\mathcal{L}(g_p)\\
               &=0.
\end{align*}
Observe that $B(0,X)=-2iIm(\bar k_1 k_2)I$ is a constant matrix.
Thus for any $t>0$,
$$
B(t, X(t, X_0))=B(0, X_0)=-2i Im(\bar k_1 k_2)I.
$$
 The condition $Im(\bar k_1 k_2)\not=0$ ensures that $g_p$ must be invertible for all $t>0$.
 Otherwise there would be a nonzero vector $c$ such that $g_p
 c=0$, hence $\bar c^\top Bc=\overline(g_p c)^\top g_xc=0$, leading to a contradiction.\\

c) Set $Q=g_x +\tilde Mg_p$. A calculation using (\ref{mp}),
(\ref{p1}) and (\ref{p2}) gives
$$
\mathcal{L}[Q]=\mathcal{L}(g_x)+\mathcal{L}(\tilde
M g_p)=-(H_{xp}+\tilde M H_{pp})Q.
$$
If $M_{\rm in}=-g_x (g_p)^{-1}$ initially, then $Q(0, X)=0$ for all
$X\in \mathbb{R}^{2n}$. Thus we have
$$
Q(t, X)=g_x + \tilde M g_p \equiv 0.
$$
This gives
$$
\tilde M=- g_x (g_p)^{-1}
$$
for all $t>0$ since $g_p$ is invertible. \\
d. i) Initially $M_{\rm in}=-\frac{k_1}{k_2}I=M_{\rm in}^T$. Since $\tilde M^T$
also satisfies equation (\ref{mp}),
 hence $\tilde M= \tilde M^T$.\\
ii) With the definition of $B$, we have
$$
B= -\overline{g_p}^T \tilde M g_p +\overline{g_p}^T
\overline{ \tilde M} g_p=-2iIm[\overline{g_p}^T \tilde M
g_p].
$$
Initially we have $g_p=k_2I$. This together with $B(t, X)=B(0,
X_0)$ along $\Gamma$ gives $Im[\overline{g_p}^T \tilde M g_p]
=|k_2|^2Im[M_{\rm in}]$. Note that $Im[M_{\rm in}]$ is positive definite, hence
$Im[\tilde M]$ remains positive definite for all $t>0$.
\end{proof}
\begin{rem}
The formula $ M=-\phi_x (\phi_p)^{-1}$,  first derived in
\cite{JLOTa:2005}, plays an important role in \cite{JLOTa:2005}
in deriving the equation
$$
\mathcal{L}[f]=0
$$
for the quantity $f(t, X)=|\tilde A(t, X)|^2 \det(\phi_p)$,
which remains globally bounded even when $\phi_p$ becomes singular.
We note that a complex level set function was used in \cite{JWY:2008} to obtain
a globally bounded Hessian.
\end{rem}
\begin{rem}
Since the Liouville equation is geometric and homogeneous, for
each fixed $X_0$, the shift $X_0$ in the level set function can be
simply ignored, and be added back whenever it is needed. In other
words, we can take initial data  $\phi(0, X)=X$, then the curve
$\Gamma$ can be represented as $X_0$ level set:
$$
\Gamma=\{X, \quad \phi(t, X)=X_0\}.
$$
\end{rem}

\begin{rem}If we follow this construction, the initial data for
$\tilde M_{\rm in}$ then depends on how we initialize the level set function
$\phi$. If we take $k_1=\beta>0$ and $k_2=i$, then $g =\beta \phi_1 +i\phi_2$.
If $(\phi_1, \phi_2)(x, X)=(x, p)$, then $\tilde
M_{\rm in}=i\beta I$. If $p_0$ is restricted to be the phase gradient at $x_0$
initially, then  the initial level set function can be chosen as
$(\phi_1, \phi_2)(0, X)=(x, p-\nabla_x S_{\rm in}(x))$, this leads to $
\tilde M_{\rm in}=\partial_x^2 S_{\rm in}(x) +i\beta I.$ In this case $\phi_2$ is the function $w(t, X)$
from (\ref{w}).
\end{rem}

\subsection{Evaluation of higher order derivatives of the phase}
Let $G(t, y)=\partial_t \Phi +H(y, \nabla_y \Phi)$. If one wants to have $G(t, y)$ vanish to a higher order than two
on $\gamma$, it is necessary to obtain higher order derivatives of
$\Phi$. In the appendix we derive a system of linear equations
for $m_\alpha(t, X)=\partial_y^\alpha \Phi(t, x(t, X_0))$ on $\gamma$.
We now show that this again can be done through the vector-valued level set function $\phi$.

Differentiating $\phi_l (t, x, \nabla_x \Phi)=0$, $l=1,2$, to order of $r\geq 3$ we obtain
$$
\sum_{j=1}^n \partial_{p_j}\phi_l \partial_{y_j}(\partial_y^\alpha
\Phi) +\sum_{|\beta|=r} c_{l \alpha \beta}\partial_y^\beta \Phi
+d_{l \alpha}=0
$$
for all multi-indices $\alpha$ of length $r$.  Let
$g=k_1\phi_1+k_2\phi_2$, again using the invertibility of
$g_p$ we can obtain
$$
\nabla_x (m_\alpha)=-(g_p)^{-1} \left[
\sum_{|\eta|=r} c_{\alpha \eta}m_\eta(t, X) + d_\alpha
\right].
$$
We do this recursively, since the coefficients $c_{\alpha,
\eta}=(c_{1\alpha, \eta}, \cdots,c_{n\alpha, \eta})^\top $ and
$d_\alpha=(d_{1\alpha}, \cdots d_{n\alpha})^\top $ depend on all the
partials up to order $r-1$. Since $g_p$ is invertible, the
obtained derivatives remain bounded for all $t>0$.

\section{Appendix }
In this appendix we follow \cite{Ralston05} to determine higher order derivatives of phase and amplitude
on $\gamma$, and further derive phase space equations they satisfy.  From $\partial_y^\alpha G=0$
on $\gamma$ with $|\alpha|\geq 3$, we obtain
$$
\partial_t (\partial_y^\alpha \Phi) +H_p\cdot \nabla_x(\partial_y^\alpha \Phi)+\sum_{|\eta|=|\alpha|}
c_{\alpha, \eta}\partial_x^\eta \Phi+ d_\alpha=0,
$$
where $c_{\alpha, \eta}$ and $d_\alpha$ depends on $\partial_y^\kappa \Phi$ for $|\kappa|<|\alpha|$.
Using the Hamiltonian equations $\frac{d}{dt}x =H_p$ we obtain
$$
\frac{d}{dt} (\partial_y^\alpha \Phi(t, x(t, X_0); X_0)+\sum_{|\eta|=|\alpha|}
c_{\alpha, \eta}\partial_x^\eta \Phi(t, x(t, X_0); X_0)+ d_\alpha=0
$$
on $(t, x(t, X_0))$.  Following Lemma \ref{lem2.1} we obtain a linear system of Liouville type PDEs for partial derivatives  $m_\alpha(t, X)$ of a fixed order:
\begin{equation}\label{liu+}
    \mathcal{L}[m_\alpha]  +\sum_{|\eta|=|\alpha|}
c_{\alpha, \eta}m_\eta + d_\alpha=0.
\end{equation}
We solve the system (\ref{liu+}) starting from $|\alpha|=3$, then $|\alpha|=4$ and so on until $|\alpha|=k+1$ for $k^{th}$ order GB solutions. Since equations are linear, we have solutions defined for all $t>0$. This construction ensures that
\begin{equation}\label{gt}
G(t, y)=O(|y-x|^{k+2}).
\end{equation}

To determine the Taylor series of $A_l, l=0, \cdots, N$ on $\gamma$, one proceeds as follows.
Define the coefficients $c_j(t, y)$ by
\begin{equation}\label{pap}
   P[A(t, y)e^{i\Phi(t, y)/\ep}]=\left( \sum_{j=0}^{N+2} c_j(t, y)\ep^j\right)e^{i\Phi/\ep}, \quad A = \sum_{j=0}^N
A_j\ep^{j}.
\end{equation}
Then, with $P=-i\ep \partial_t +H(y, -i\ep \partial_y)$, we obtain
\begin{align*}
c_0(t, y) &= G(t, y)A_0,\\
c_1(t, y) &=-iLA_0 +G(t, y)A_1\\
c_{l+1}(t, y) &=-iL A_{l}+G(t, y)A_{l+1}+g_l, \quad l=1, \cdots N+1,
\end{align*}
where $L$  is a linear differential operator with
coefficients depending on $\Phi$,
$$
L= \partial_t+ H_p\cdot \nabla_y +
\frac{1}{2}\left[ tr(H_{yp})+tr(M H_{pp}(y, \Phi_y))\right],
$$
and $g_l=-\frac{1}{2}\Delta_y A_{l-1}$.

Thus to make $P[\psi^\ep] = O(\ep^{K})$ for a given $K\in \mathbb{Z}$, we now only need to make $c_j$ vanish on $\gamma$ to
sufficiently high order.  To do so  we can solve the equations $L A_l + ig_l = 0$ recursively starting with
$l = 0$ ($g_0=0$), and solve it to arbitrarily high order by solving the linear
transport equations for the partial derivatives of $A_l$ that one
gets by differentiating the above equations.   From the above procedure we see that the number of terms, $N$, in the solution ansatz for $k^{th}$order GB approximation is determined by the following relation
$$
\frac{k-1}{2}<N \leq \frac{k+1}{2}
$$
In other words, $N=\lfloor\frac{k-1}{2}\rfloor+1$. Actually, given $G(t, y)$ vanishes to order $k+1$ on $\gamma$, we can choose the Taylor series of $A_0$ on $\gamma$ up to order $k-1$ so that $c_1$ vanishes to order $k-1$ on $\gamma$. Passing to the higher order equations $c_{l+1}=0$, we see that we can choose $A_l$ so that $c_{l+1}$ vanishes on $\gamma$ to order $k+1-2(l+1)$. Thus we need $2>k+1-2N \geq 0$.

Thus for $k^{th}$ order GB solutions, it is necessary to compute $\partial_y^\alpha A_l$ for $|\alpha|\leq k-1-2l$:
{\small
\begin{align*}
& L(\partial_y^\alpha A_0)+\sum_{|\eta|<|\alpha|} \left(
                                                         \begin{array}{c}
                                                           \alpha \\
                                                           \eta \\
                                                         \end{array}
                                                       \right)
 \partial_y^{\alpha-\eta}L (\partial_y^\eta A_0)\Big|_{\gamma}=0, \quad |\alpha|\leq k-1, \\
& L(\partial_y^\alpha A_l)+\sum_{|\eta|<|\alpha|} \left(
                                                         \begin{array}{c}
                                                           \alpha \\
                                                           \eta \\
                                                         \end{array}
                                                       \right)
 \partial_y^{\alpha-\eta}L (\partial_y^\eta A_l) -\frac{i}{2}\Delta_y  \partial_y^\alpha A_{l-1}
 \Big|_{\gamma}=0, \quad |\alpha|\leq k-1-2l, \; l=1, \cdots \left \lfloor \frac{k-1}{2} \right\rfloor.
\end{align*}
}
Lifting the operator into the phase space, we can obtain  $\tilde A_l(t, X)$  recursively by
solving
$$
\mathcal{L}[\tilde A_l]=- \frac{\tilde A_l}{2}\left[ tr(H_{yp})+tr(\tilde M(t, X) H_{pp}(y, p))\right] - \tilde g_l,
$$
where $\mathcal{L}$ is the Liouville operator.  Same lifting can be applied to all involved derivatives of the amplitude $A_l$ for $l=0, \cdots \left \lfloor \frac{k-1}{2} \right\rfloor.$ This completes the construction for all involved Taylor coefficients of both phase and amplitude.

\bigskip
\section*{Acknowledgments} H. Liu wants to thank the Department of Mathematics at UCLA for its hospitality and support during his visit in winter quarter of 2009 when this work was completed.  Liu's research was partially supported by the National Science Foundation under the Kinetic FRG Grant DMS07-57227.

\bigskip

\end{document}